\documentclass[a4paper,12pt]{article}
\usepackage{amsmath}
\usepackage{amsfonts}
\usepackage{amsthm}
\usepackage{amssymb}
\usepackage{yfonts}
\usepackage{[leqno}
\title{On a Generalization of the van der Waerden Theorem}
\author{Rudi Hirschfeld\\
University of Antwerp\\
\texttt{rudihirschfeld@hetnet.nl}}
\newtheorem{theorem}{Theorem}
\newtheorem{corollary}{Corollary}

\newtheorem{bonus}[theorem]{Bonus}

\newcommand{\betaN}{\beta\mathbb{N}}

\newcommand{\N}{\mathbb{N}}
\newcommand{\betaomega}{\beta\omega}
\newcommand{\SC}{Stone-\v{C}ech \mbox{}}

\usepackage{natbib}
\numberwithin{equation}{section}
\begin{document}
\maketitle

\begin{abstract}
For a given length and a given degree and an arbitrary partition of the positive integers, there always is a cell containing a polynomial progression of that length and that degree; moreover, the coefficients of the generating polynomial can be chosen from a given subsemigroup and one can prescribe the occurring powers. A multidimensional version is included.
\end{abstract}
\section{Introduction}
A sequence in $\mathbb{R}$ will be called a \emph{polynomial progression} if it is of the form $\{P(1),P(2),P(3),\ldots\}$ for some polynomial $P(x)=a_dx^d+a_{d-1}x^{d-1}+\cdots +a_1x+a_0$. This progression is said to be of \emph{degree} $d$ if $P$ has degree equal to $d$ and not less.

\begin{theorem}
Given two positive integers $d\mbox{ and }l$, if the set of the positive integers is split up into finitely many non-overlapping parts, there exists a polynomial progression of length $l$ and of  degree $d$ that belongs to precisely one of these parts.
\end{theorem}

For $d=1$ the polynomials look like $P(x)=a+bx$ and the $l$-segment of the polynomial progression takes the form $\{a+b,a+2b,a+3b,\ldots,a+lb\}$: the theorem boils down to the well-known van der Waerden Theorem on monochromatic arithmetic progressions. It is fun to write down the $d=2$ case.
\begin{corollary}
Given any $l\in \N$ and any finite coloring of $\N$, there exist three positive integers $a,b$ and $c$ for which all terms in $\{a+b+c,a+2b+4c,a+3b+9c,\ldots,a+lb+l^2c\}$ have the same color.
\end{corollary}

The 1927 proof of van der Waerden's Theorem is quite complicated, involving a double induction argument. The 1927 issue of the journal, \cite{vdW}, is difficult to access nowadays, but a very clear exposition is found in \textsc{R.L. Graham, B.L. Rothschild and J.H. Spencer} \cite{GraRotSpe}, pp 29 -- 34. As B.L. van der Waerden once remarked, around 1927 he was not aware of the impact of his result as a prototypical Ramsey Theorem - after all, Ramsey's famous paper stems from 1930 - and merely considered it as a clever exercise. A proof of the above theorem by means of induction seems a Sisyphean task. We rather use some ideal theory in the semigroup $\betaN$. As a matter of fact, the argument in the \textsc{Hindman-Strauss} treatise \cite{HinStr} for the van der Waerden Theorem (see 14.1 \emph{l.c.}) is readily adapted to the present situation. By preferring the smooth $\betaN$-argument to a complicated induction proof we ignore the calvinistic concern (see  \cite{HinStr} p.280) that it ``is enough to make someone raised on the work ethic feel guilty''.

The more restrictions one puts on the admissible polynomials, the fewer polynomials one has at his/her disposal and the more difficult it seems to force the ensuing polynomial progressions into one and the same cell. The polynomials $P(x)=\sum_{k=0}^d a_kx^k$ we admit here satisfy
\begin{itemize}
	\item the admissible coefficients $a_k$ belong to one and the same subsemigroup $\mathbb{S}$ of $(\omega,+)$, where $\omega=\N \cup {0}$;
	\item the admissible exponents in the powers $x^k$ belong to a subset $\mathbb{D} \subset \{0,1,2,\ldots,d\}$ containing $d$.
\end{itemize}
Such polynomials will be called $(\mathbb{S},\mathbb{D})$-polynomials.

The sharpened theorem reads

\begin{theorem}
Given two positive integers $d\mbox{ and }l$, if the set of the positive integers is split up into finitely many non-overlapping parts, there exists a polynomial progression of length $l$ and of  degree $d$, generated by a $(\mathbb{S},\mathbb{D})$-polynomial, that belongs to precisely one of these parts.
\end{theorem}

\section{Proof}
Since Theorem 1 concerns the special case where $\mathbb{S}=\omega$ and $\mathbb{D}=\{0,1,2,\ldots,d\}$, we only need to prove Theorem 2.

Fix $d$ and $l$ in $\N=\{1,2,3\ldots\}$. Without los of generality we may assume that $l>d$. In fact, once the theorem has been proved for ``long'' progressions (that is $l>d$), then the pertinent cell certainly contains shorter segments ($l\leq d$). We consider polynomials $P(x)=\sum_{i=0}^d a_ix^i$ in one indeterminate $x$ of degree $\leq d$ with coefficients in $\omega^{d+1}$.
Consider the following sets $S_o$ and $I_o$ in $\omega^l$ consisting of $l$ consecutive polynomial values                                                                                                                                                                                                                                                                                                                                                                                                                                                                                                                                                                                                         
\begin{eqnarray*}
S_o &=&\{\{P(1),P(2),\ldots,P(l)\}\in\omega^{l}:\;P(x)=\sum_{k\in \mathbb{D}} a_k x^k, \mbox{ with }\{a_0,a_1,\ldots,a_d\} \in \mathbb{S}^{d+1}\}\\
I_o &=&\{\{P(1),P(2),\ldots,P(l)\}\in\N^{l}:\;P(x)=\sum_{k\in \mathbb{D}}a_k x^k, \mbox{ with }\{a_0,a_1,\ldots,a_d\} \in (\mathbb{S}\cap\N)^{d+1}\}\\
\end{eqnarray*}
The impact of the assumption that $l>d$ is that each element in $S_o$ corresponds to a \emph{unique} polynomial. In fact, if such an $l$-tuple would be generated by two different polynomials, the difference of these polynomials would have more zeros (\emph{viz.} at the $l$ points $1,2,\ldots,l$ in $\mathbb{C}$) than its degree $d<l$ permits.

$S_o$ is a subsemigroup of $\mathbb{S}^{d+1}$ under coordinatewise addition, the restrictions $k\in\mathbb{D}$ meaning that only addition of coordinates $k$ from $\mathbb{D}$ matters. In fact, the sum of two $l$-tuples in $S_o$ corresponds to the sum of their unique polynomials and the latter is again a polynomial of degree $\leq d$ with coefficients in the semigroup $\mathbb{S}$. 

The progressions $\{P(1),P(2),\ldots,P(l)\}$ in $I_o$ all have degree $=d$, since $a_d\geq 1$. It follows that $I_o$ is a proper subset of $S_o$. Obviously, $I_o$ is also a semigroup. Moreover, $I_o$ is a ideal in $S_o$. In fact, upon adding any point in $S_o$ to an arbitrary element of $I_o$,  all coefficients of the sum polynomial are again $\geq 1$ and this polynomial is of exact degree $d$. Although trivial, we notice that $S_o$ contains constant $\N$-valued polynomials, but $I_o$ contains none of these. This will be instrumental shortly.

Consider the \SC compactification $\betaomega$. We are going to use a few facts about $\betaomega$ that are found in \textsc{N. Hindman and D. Strauss} \cite{HinStr}. We find it convenient to ignore the slight differences in the ideal theory between the two  \emph{semigroups} (see \cite{HinStr}, Chap. 4) $\betaomega$ and $\betaN$, writing $\betaN$ where $\betaomega$ would sometimes be more appropriate. From this point on we can follow the proof of the van der Waerden Theorem in \cite{HinStr}, Theorem 14.1, almost \emph{verbatim}. 

Take the compact product space $Y=(\betaN)^l$ and the closures $S=cl_Y(S_o)$ and $I=cl_Y(I_o)$. The semigroup $\betaN$ has a \emph{smallest} ideal $K(\betaN)\neq \emptyset$ (see \cite{HinStr},  Chap 4), which will be our main tool.

Take any point $p \in K(\betaN)$ and consider the constant $l$-tuple $\vec{p}=\{p,p,\ldots,p\}$. The crucial step is to show that $\vec{p}$ belongs to $S$. 

The closures $cl_{\betaN}B$ of the members $B\in p$ form a neighborhood basis in $\betaN$ around $p$. It follows that for the product topology in $Y$ there exist members $B_1,B_2,\ldots,B_r \in p$ for which the box $U=\prod_{1\leq i \leq r}cl_{\betaN}(B_i)$ is a $Y$-neighborhood of $\vec{p}$. The intersection $\cap_{1\leq i \leq r}cl_{\betaN}(B_i)$ is a $\betaN$-neighborhood of $p$. The set $\N$ lying dense in $\betaN$, it is intersected by this neighborhood. Select $a \in \N \cap \big(\cap_{1\leq i \leq r}cl_{\betaN}(B_i)\big)$. The constant $l$-string $\vec{a}=\{a,a,\ldots,a\}$ thus belongs to $U$. Also, $S_o$ containing all constant $l$-tuples, we have $\vec{a} \in S_o$. Consequently, we have $\vec{a}\in S_o \cap U$. This shows that $\vec{p}$ belongs to the closure of $S_o$ in $Y$, and so $\vec{p} \in S$, indeed. 

Next we use the fact that by \cite{HinStr}, Theorem 2.23,  the $K$-functor preserves products. From $p\in K(\betaN)$ we infer $\vec{p} \in \big(K(\betaN)\big)^l=K\big((\betaN)^l\big)=K(Y)$. Conclusion: $\vec{p}\in S \cap K(Y)$.

Having shown that $S \cap K(Y) \neq \emptyset$, we can invoke \cite{HinStr}, Theorem 1.65 to determine the smallest ideal of the semigroup $S$: it simply is $K(S)=S \cap K(Y)$.
This leads to 
\begin{equation}
\vec{p} \in K(S). \label{bijna klaar}
\end{equation}

Obviously, $I$ is an ideal in $S$. The smallest ideal in $S$ is contained in $I$: $K(S)\subset I$. It follows from (\ref{bijna klaar}) that $\vec{p}\in I$.

Finally, let $\N=\bigcup_i A_i$ be a finite partition. The closures $\bar{A}_i=cl_{\betaN}A_i$ are open and form a partition of $\betaN$. Hence, precisely one of them, $\bar{A}_j$ say, is a $\betaN$-neighborhood of our point $p \in K(\betaN)$. Then $V=\big(\bar{A}_j\big)^l$ is a $Y$-neighborhood of $\vec{p}$. Because $\vec{p}\in I$, $V$ must meet the dense subset $I_o$ of $I$ and we can select a polynomial $P$ in in such a manner that $\{P(1),P(2),\ldots,P(l)\}$ belongs to $V$. But the $P(1),P(2),\ldots,P(l)$ still are integers in $\N$. For this reason
\[\{P(1),P(2),\ldots,P(l)\}\subset \bar{A}_j \cup \N=A_j\]
and the segment $\{P(1),P(2),\ldots,P(l)\}$ has the color of $A_j$. \qed
 \vspace{6mm}

\section{Free gifts}

The essential property of the set $A_j$ used in the last part of the above proof is the fact that $\bar{A}_j$ contains a point $p$ belonging to $K(\betaN)$, or $A_j \in p$. Sets $A \subset \N$ belonging to some $p\in K(\betaN)$ are called \emph{piecewise syndetic} sets. We recall that in terms of $\N$ itself, $A$ is piecewise syndetic if and only if the \emph{gaps} between its intervals of consecutive elements remain bounded in lengths, (see \cite{HinStr} Theorem 4.40). It follows that $A_j$ may be replaced by any infinite piecewise syndetic set $A$ and we get as a

\begin{bonus} Given a piecewise syndetic set $A\subset\N$, a length $l$ and a degree $d$, there exists a polynomial progression of degree $d$ for which the first $l$ terms belong to $A$.
\end{bonus}

Finally we consider a \emph{multidimensional} version of the theorem, dealing with $m$ polynomial progressions of varying lengths and degrees simultaneously.
\begin{bonus}
Pick the following items in $\N$: a dimension parameter $m$, degrees $d_1,d_2,\ldots,d_m$, and lengths $l_1,l_2,\ldots,l_m$ . If the set $\N$ is split up into finitely many non-overlapping parts, there exist $m$ polynomial progressions of length $l_i$ and of  degree $d_i$ each, $1\leq i \leq m$, that simultaneously belong to  one of these parts. Also, any given piecewise syndetic set contains such a collection of polynomial progressions.
\end{bonus}

\noindent \textbf{Remark} There is an obvious $(\mathbb{S},\mathbb{D})$ version.

\noindent\emph{Proof.}
We introduce arrays 
\begin{gather*}
\mathcal{P}=
\begin{pmatrix}P_1(1) & P_1(2) \cdots P_1(l_1)\\
P_2(1) & P_2(2) \cdots P_2(l_2)\\
\cdots &\cdots \cdots \cdots\\
P_m(1) & P_m(2) \cdots P_m(l_m)
\end{pmatrix}
\end{gather*}
generated by polynomials $P_1,P_2,\ldots,P_m$ with coefficients from $\omega$. 

These arrays $\mathcal{P}$ need not have the customary rectangular form, the $i^{th}$ row having $l_i$ entries. Extending these rows by putting zeros in the empty places until they all get $\max\{l_i: i=1,2,\ldots,m\}$ entries would unnecessarily complicate the definition of $\mathcal{I}_o$ \emph{infra}.) 

We have avoided to call these $\mathcal{P}$ matrices since they are not intended to act as transformations in some vector space. In order to describe the set these arrays belong to we write $e_i=\{0,\ldots,1,0,\ldots,0\}=\delta_{ij}$ for the usual unit vectors in $\mathbb{R}^m$. These unit vectors are customarily envisaged as rows; upon transposition we get the unit columns $e_i^T$ . The direct sum decomposition 
\[\omega^m=e_1^T \omega\oplus e_2^T\omega\oplus \cdots \oplus e_m^T \omega\]
divides $\omega^m$, and thereby $\N^m$, into $m$ horizontal layers, each equal to $\N$ and each row is an additive semigroup at its own. 
\begin{gather*}
e_i ^T \N=
 \begin{pmatrix}
 	0&\\
	\;\;\cdots\\
	1&\hspace{-3mm}2\;\;3\;\;4\;\;5\;\;\ldots\\
	0&\\
	\;\;\cdots\\
	0&
\end{pmatrix}\\
\end{gather*}
Picture: the $i^{th}$ row of an array $\mathcal{P}$ is contained in the $i^{th}$ layer. 

Upon replacing the $l$-tuples in the definitions of $S_o$ and $I_o$ by the arrays $\mathcal{P}$ we get 
\begin{eqnarray*}
\mathcal{S}_o &=&\{\mathcal{P}\in\bigoplus_{i=1}^m e_i^T \N:\; P_i(x)=\sum_{k=0}^{d_i} a_{ki} x^k, \mbox{ with }\{a_0,a_1,\ldots,a_{d_{i}}\} \in \omega^{d_{i}+1} \mbox{ for }1=1,2,\ldots,m\},\\
\mathcal{I}_o &=&\{\mathcal{P}\in \bigoplus_{i=1}^m e_i^T \N:\; P_i(x)=\sum_{k=0}^{d_i} a_{ki} x^k, \mbox{ with }\{a_0,a_1,\ldots,a_{d_{i}}\} \in \N^{d_{i}+1} \mbox{ for }1=1,2,\ldots,m\}.\\
\end{eqnarray*}
These are subsemigroups of the $\bigoplus_{i=1}^m e_i^T \N$ and $\mathcal{I}_o$ is a proper ideal in $\mathcal{S}_o$. 

We refrain from repeating all details the above proof for the $m=1$ case. 

For a start, we may assume without loss of generality that $l_1>d_1,l_2>d_2,\ldots,l_m>d_m$. Define $l=\max_{1\leq i \leq m}l_i$. This time we have to deal with the compact space $Y^m$, one $Y=(\betaN)^l$ for each layer, so that $Y^m=(\betaN)^{lm}$. The closure 
$\mathcal{I}=cl_{Y^{m}}(\mathcal{I}_o)$ is an ideal in the semigroup $\mathcal{S}=cl_{Y^{m}}(\mathcal{S}_o)$

To every $p\in K(\betaN)$ we assign the constant $m\times l$ array
\begin{gather*}
\vec{\mathbf{p}}=
\begin{pmatrix}p & p\cdots p\\
p & p \cdots p\\
&\hspace{-6mm}\cdots \cdots \cdots \cdots\\
p & p\cdots p
\end{pmatrix}
\end{gather*}
After a little twist the above argument leads to $\vec{\mathbf{p}} \in K(\mathcal{S})\subset \mathcal{I}$. For a piecewise syndetic set $A\in p$ the product $V=(\bar{A})^l$ is a $Y^m$-neighborhood of $\vec{\mathbf{p}}$ which intersects the dense subset $\mathcal{I}_o$ of $\mathcal{I}$ in at least one point. This point is an array $\mathcal{P}$, say. It follows that the entries $P_i(j)$ of $\mathcal{P}$ belong to $cl_{Y^m}A$ and thus to $cl_{\betaN}A$.  All $P_i(j)$ being positive integers, we may write
\[\bigcup \{P_i(j):i=1,\ldots,m;j=1,\ldots l_i\} \subset \bar{A} \cap \N.\]
Conclusion: these $m$ polynomial progressions do lie in $A$ itself. \qed

\end{document}